\documentclass[10pt]{article}
\usepackage{amsmath,amsthm,amsfonts}

\input epsf.tex
\newdimen\epsfxsize
\newdimen\epsfysize

\theoremstyle{plain}
\newtheorem{thm}{Theorem}[section]

\newtheorem{conj}{Conjecture}[section]

\begin{document}

\title{Lower Bounds on Virtual Crossing Number and Minimal Surface Genus}

\author{Kumud Bhandari \\ H. A.  Dye \\
Louis H. Kauffman   }

\maketitle

\begin{abstract}
We compute lower bounds on the virtual crossing number and minimal surface genus of virtual knot diagrams from the arrow polynomial. In particular, we focus on several interesting examples. 
\end{abstract}

\section{The arrow polynomial}

The arrow polynomial is an invariant of oriented link diagrams introduced in \cite{dk-arrow}. This polynomial is an element of the ring $ \mathbb{Z} [[ A, A^{-1}, K_1, K_2, \ldots ]] $ and is invariant under both the classical and virtual Reidemeister moves. This polynomial is equivalent to the simplified extended bracket \cite{louskein} and the Miyazawa polynomial \cite{miya}. 

\begin{figure}[htb] \epsfysize = 1.25 in
\centerline{\epsffile{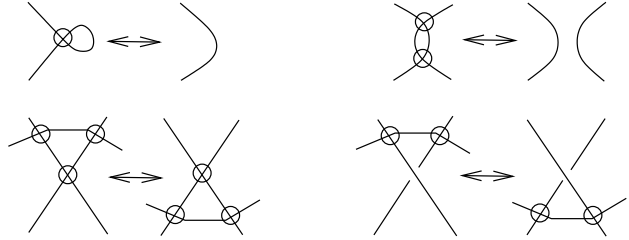}}
\caption{VirtualReidemeister moves}
\label{fig:vrmoves}
\end{figure}

A \textit{virtual link diagram} is a decorated immersion of $n$ copies of $S^1 $ into the plane, with two types of double points: virtual and classical crossings. Classical crossings are indicated by over under markings and virtual crossings are indicated by a solid, circled crossing.
A \textit{virtual link} is an an equivalence class of virtual link diagrams. Two virtual link diagrams are equivalent if one can be transformed into the other by a sequence of classical and virtual Reidemeister moves. 
\begin{figure}[htb] \epsfysize = 1.25 in
\centerline{\epsffile{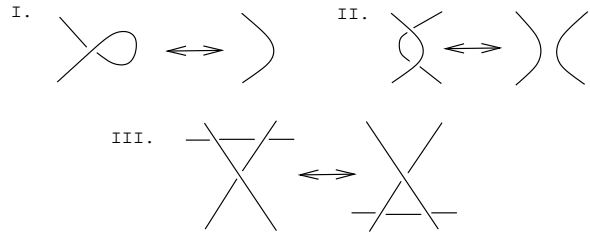}}
\caption{Classical Reidemeister moves}
\label{fig:rmoves}
\end{figure}
The classical Reidemeister moves are shown in figure \ref{fig:rmoves} and the virtual Reidemeister moves are shown in figure \ref{fig:vrmoves}.

\begin{figure}[htb] \epsfysize = 0.7 in
\centerline{\epsffile{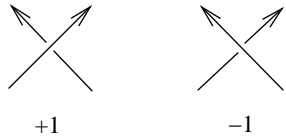}}
\caption{Crossing sign }
\label{fig:sgn}
\end{figure}

An oriented virtual link diagram is determined by assigning an orientation to each component of the diagram. For an $n$ component link, there are $2^n$ possible orientations. Based on the orientation of each component, we can determine a numerical value associated with each classical crossing $v$ as shown in figure \ref{fig:sgn}. This numerical value associated with a classical crossing $v$ is called the \textit{crossing sign} and denoted $sgn(v)$. Based on the crossing sign, we can compute the \textit{writhe} of a link diagram $L$, denoted $w(L)$. The writhe is computed by summing over all classical crossings $v$ in the diagram. That is,
\begin{equation}
w(L) = \sum_{v} sgn(v)
\end{equation}.

\begin{figure}[htb] \epsfysize = 2.0 in
\centerline{\epsffile{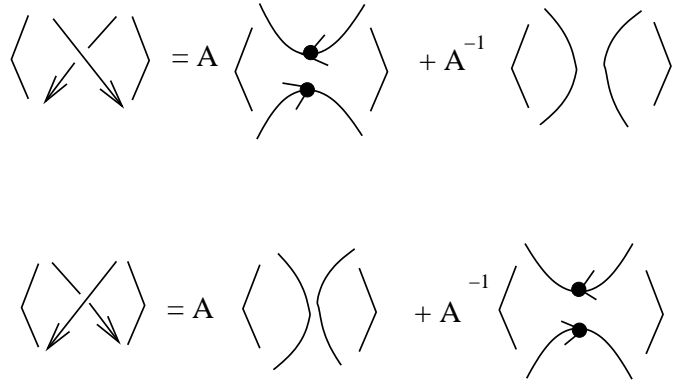}}
\caption{Arrow polynomial skein relation}
\label{fig:arrowstate}
\end{figure}

The arrow polynomial, defined in \cite{dk-arrow}, is an invariant of oriented, virtual link diagrams. These polynomials are elements of the commutative ring:
\begin{equation*}
 \mathbb{Z} [[ A, A^{-1}, K_1, K_2, \ldots ]]
\end{equation*}
where the $K_i$ form an infinite set of variables. The arrow polynomial is obtained from the oriented skein relation shown in figure \ref{fig:arrowstate}.
\begin{figure}[htb] \epsfysize = 1.0 in
\centerline{\epsffile{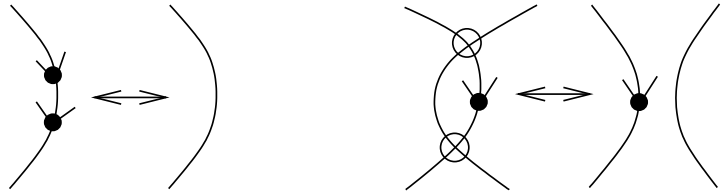}}
\caption{Virtual Reidemeister move with a nodal arrow}
\label{fig:reduce}
\end{figure}
\begin{figure}[htb] \epsfysize = 1.0 in
\centerline{\epsffile{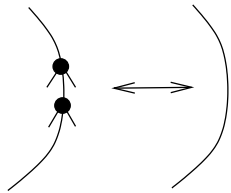}}
\caption{Reducing the total number of arrows}
\label{fig:reduce2}
\end{figure}
Applying this skein relation to each classical crossing results in a weighted sum of states with coefficients in $ \mathbb{Z} [[A, A^{-1}]]$.We obtain a state of the diagram by choosing either a horizontal or vertical smoothing for each classical crossing in the diagram. This state consists of a collection of closed loops, possibly with virtual crossings. Each loop contains a non-negative, even number of nodal arrows. We can simplify each state into a collection of disjoint loops possibly containing nodal arrows by applying the virtual Reidemeister moves to the diagram and using the move shown in figure \ref{fig:reduce}.
The total number of nodal arrows in a component is reduced using the convention shown in figure \ref{fig:reduce2}. 
We evaluate a single loop $C$ with $2n$ arrows as follows:
\begin{equation}
 \langle C \rangle = K_n.
 \end{equation}
 
\begin{figure}[htb] \epsfysize = 1.5 in
\centerline{\epsffile{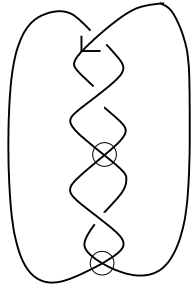}}
\caption{Virtualized trefoil}
\label{fig:vtrefoil}
\end{figure}

Then for a state $S$ with 
\begin{equation*}
S= \prod_{i=1} ^{n} C_i 
\end{equation*}
we find that 
\begin{equation}
\langle S \rangle = \prod_{i=1} ^{n} \langle C_i \rangle.
\end{equation}
Let $d=-A^2 -A^{-2}$ and let $L$ denote a virtual link, then the arrow polynomial of $L$ is:
\begin{equation}
\langle L \rangle_A = \sum_{S}  A^{ \alpha - \beta} d^{|S|-1} \langle S \rangle
\end{equation}
where $ \alpha $ is the number of smoothings in the state $S$ with coefficient $A$ and $ \beta $ is the number of smoothings with coefficient $ A^{-1}$, and $|S| $ denote the number of closed loops in the state. Recall that$ \langle L \rangle_{A} $ is invariant under the virtual Reidemeister moves and the classical Reidemeister moves II and III \cite{dk-arrow}. 
The normalized arrow polynomial, denoted $ \langle L \rangle_{NA} $, is invariant under all classical and virtual Reidemeister moves. The normalized arrow polynomial of a link $L$ is:
\begin{equation}
(-A^3)^{- w(L)} \langle L \rangle_A.
\end{equation}
For example, the arrow polynomial of the knot shown in figure \ref{fig:vtrefoil} is:
\begin{equation}
-A^{-3} ( -A^{-5} + K_1 ^2 A^{-5} - K_1^2 A^3 ).
\end{equation}

The arrow polynomial determines a lower bound on both the genus and the virtual crossing number of a virtual link.
Recall that the \textit{ virtual crossing number} of a link $L$ (denoted $v(L)$) is the minimum number of virtual crossings in any virtual link diagram equivalent to $L$. Notice that individual summands of $ \langle L \rangle_A $ have the form:
\begin{equation*}
A^m K_{i_1} ^{p_1} K_{i_2} ^{p_2}  \ldots K_{i_n} ^{p_n}
\end{equation*}
The \textit{k-degree} of the summand is:
\begin{equation*}
i_1  \times p_1 + i_2 \times p_2 + \ldots + i_n \times p_n
\end{equation*}
The \textit{maximum k-degree} of $ \langle L \rangle_A $ is the maximum k-degree of the summands. 
In \cite{dk-arrow}, the following theorem was proved.
\begin{thm}Let $K$ be a virtual link diagram. Then the virtual crossing number of $K$, $ v(K) $, is greater than or equal to the maximum k-degree of $ \langle K \rangle_A $. 
\end{thm}
Hence the maximum k-degree provides a lower bound on the virtual crossing number.

Recall that virtual links are in one to one correspondence with representations of virtual links ( \cite{kamada-stable}, \cite{dk-surf}). A representation of  a virtual link, denoted $(F,L)$, is an embedding of the link $L$ in $ F \times I $ where $F$ is a closed, two dimensional, oriented surface modulo Dehn twists, isotopy of the link with in $F \times I $, and handle addition/subtraction. Representations are in one to one correspondence with virtual links. 
Kuperberg proved the following theorem in \cite{kup}:
\begin{thm}Every stable equivalence class of links in thickened surfaces has
a unique irreducible representative. \end{thm}

Hence, each virtual link corresponds to a representation with a minimum genus surface. In \cite{dk-arrow}, it was shown that the arrow polynomial can determine a lower bound on the mimimum genus.
\begin{thm}\label{g} Let $S$ be an oriented, closed, 2-dimensional surface with genus $g \geq 1$. If $g =1$ then $S$ contains at most $1$ nonintersecting, essential curve and if $ g > 1 $ then $S$ contains at most $3g-3$ non-intersecting, essential curves.
\end{thm}
That is, if $ \langle L \rangle_A $ contains a $K_i$ in any summand, the minimum genus is at least one. If a single summand contains a factor of the form $K_i K_j$ then the minimum genus is at least two. 

\section{Computations}

In the following tables, the arrow polynomial has been computed for all virtual knots with four classical crossings. This tabulation is based on Jeremy Green's knot tables, available at:
http://www.math.toronto.edu/drorbn/Students/GreenJ/   \cite{jg}. Images of all the knots in this paper are available at this website. Since virtual knots have only one component, it is not necessary to specify the orientation of the links. From the arrow polynomial, we have computed both a lower bound on the virtual crossing number and the genus which is also listed in the table. There are only four knots with arrow polynomial one: (4.46, 4.72, 4.98, 4.107). In comparison, 24 four crossing knots (out of 108 knots) have Jones polynomial one. 

The maximum lower bound on virtual crossing number is three, and based on the computational results we make the following conjecture:
\begin{conj}Given a virtual knot, $K$, an upper bound on the number of virtual crossings is determined by the minimum number of classical crossings.
\end{conj}
Unlike the virtual crossing number, the classical crossing number can be determined from the Gauss diagram.
We focus on several specific examples in the remainder of the paper. 

\subsection{Knot 4.01}

\begin{figure}[htb] \epsfysize = 1.0 in
\centerline{\epsffile{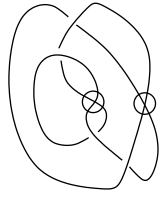}}
\caption{Knot 4.01}
\label{fig:knot401}
\end{figure}
The arrow polynomial of knot 4.01 is:
\begin{equation}
A^8 K_1^2 -3 K_1^2 +2 -2A^4 K_1^2 +2 K_2 ^1 -2A^2 K_1 +2A^{-2} K_1 + A^{-4}.
\end{equation}
The lower bound on the virtual crossing number is two and the lower bound on the genus is one. However, this minimal genus of this
virtual knot is two.

\subsection{Knot 4.09}

\begin{figure}[htb] \epsfysize = 1.0 in
\centerline{\epsffile{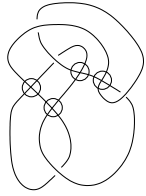}}
\caption{Knot 4.09}
\label{fig:knot409}
\end{figure}

The arrow polynomial of the knot 4.09 is:
\begin{equation*}
-A^{4} -A^{2}K_1 +A^{-2}K_1 +A^{-4}.
\end{equation*}
The lower bound on the virtual crossing number and the mimimal genus is one.
Note a single detour move reduces the number of virtual crossings to three. 
\subsection{Knot 4.22}

\begin{figure}[htb] \epsfysize = 1.0 in
\centerline{\epsffile{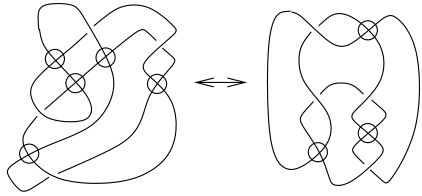}}
\caption{Knot 4.22 and an equivalent knot}
\label{fig:knot422}
\end{figure}
The arrow polynomial of the knot 4.22 is:
\begin{gather*}
-A^{6}K_1 +A^{2}K_1 +K_2 +2 -A^{2}K_1K_2 -A^{-2}K_1K_2 
      -K_1^{2} \\ -A^{-4}K_1^{2} +A^{2}K_1^{3} +2A^{-2}K_1^{3} +A^{-6}K_1^{3} -A^{-2}K_1 -A^{-6}K_1
\end{gather*}
The lower bound on the virtual crossing number is 3 and the lower bound on the minimal genus is two. A sequence of Reidemeister moves reduces the number of virtual crossings to three. The knot pictured in the right hand side of figure \ref{fig:knot422} is equivalent to knot 4.22. This demonstrates that the virtual crossing number is three.

\subsection{Knot 4.47}

\begin{figure}[htb] \epsfysize = 1.0 in
\centerline{\epsffile{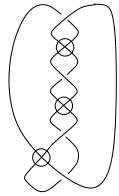}}
\caption{Knot 4.47}
\label{fig:knot447}
\end{figure}
The arrow polynomial is:
\begin{equation*}
A^{2}K_3^{1} +1 -A^{-2}K_1K_2 -A^{2}K_1K_2 +A^{-2}K_1
\end{equation*} 
The virtual crossing number is three and the minimal genus is at least two. The Jones polynomial does not differentiate between this 
knot and the unknot.

\subsection{Knot 4.91}

\begin{figure}[htb] \epsfysize = 1.0 in
\centerline{\epsffile{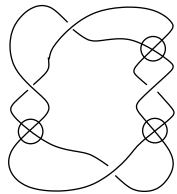}}
\caption{Knot 4.91}
\label{fig:knot491}
\end{figure}
The arrow of polynomial of knot 4.91 is:
\begin{equation*} 
 -A^{10}K_1^{3} -A^{6}K_1^{3} +A^{2}K_1^{3} +2A^{6}K_1 
       +A^{-2}K_1^{3} -2A^{2}K_1 +A^{-4}
\end{equation*}
The virtual crossing number of the knot 4.91 is three and the minimal genus is one as predicted by the arrow polynomial. 

\subsection{Knot 4.99}

\begin{figure}[htb] \epsfysize = 1.0 in
\centerline{\epsffile{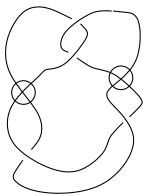}}
\caption{Knot 4.99}
\label{fig:knot499}
\end{figure}
The arrow polynomial of the knot 4.99 is:
\begin{equation*}
A^{8} -A^{4} -A^{-4} +1 +A^{-8}.
\end{equation*}
This results in a lower bound on the virtual crossing number and the minimal genus of zero. 
We observe that under virtualization, this knot is equivalent to a classical knot.

\subsection{Knots with arrow polynomial one}
The knots 4.46, 4.72, 4.98, and 4.107 have arrow polynomial one and are equivalent to the unknot via a sequence of classical and virtual Reidemeister moves and the Z-equivalence.
\begin{figure}[htb] \epsfysize = 0.5 in
\centerline{\epsffile{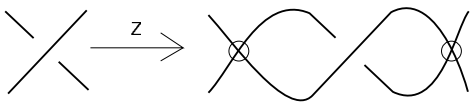}}
\caption{Z-equivalence}
\label{fig:zmove}
\end{figure}

\begin{figure}[htb] \epsfysize = 1.0 in
\centerline{\epsffile{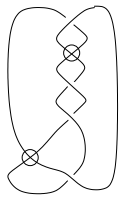}}
\caption{Knot 4.46}
\label{fig:knot446}
\end{figure}

\begin{figure}[htb] \epsfysize = 1.0 in
\centerline{\epsffile{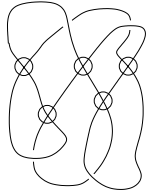}}
\caption{Knot 4.72}
\label{fig:knot472}
\end{figure}

\begin{figure}[htb] \epsfysize = 1.0 in
\centerline{\epsffile{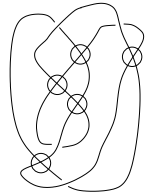}}
\caption{Knot 4.98}
\label{fig:knot498}
\end{figure}

\begin{figure}[htb] \epsfysize = 1.0 in
\centerline{\epsffile{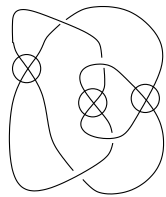}}
\caption{Knot 4.107}
\label{fig:knot4107}
\end{figure}

The knots in figures \ref{fig:knot446}, \ref{fig:knot472}, \ref{fig:knot498}, and \ref{fig:knot4107} have arrow polynomial one.
In the paper \cite{frm}, the authors (Fenn, Kauffman, and Manturov) made the following conjecture:
\begin{conj}Let $K$ be a virtual knot. If the bracket polynomial of $K$, $ \langle K \rangle =1 $ then $K$ is Z-equivalent to 
the unknot. \end{conj}
Note that if $ \langle K \rangle_A = 1 $ (arrow polynomial) then $ \langle K \rangle =1 $ (bracket polynomial). This fact and our experimental evidence support the following conjecture:
\begin{conj}
Let $K$ be a virtual knot. If $ \langle K \rangle_A =1 $ then $K$ is Z-equivalent to the unknot.
\end{conj}

\begin{table*}
		\begin{tabular}{l|p{3 in} | l | l }
			Knot & Arrow Polynomial & v(K) & g(K) \\ \hline
			4.01 & \vspace{0.1 mm} $ A^8 K_1^2 -3 K_1^2 +2 -2A^4 K_1^2 +2 K_2 ^1 -2A^2 K_1 +2A^{-2} K_1 + A^{-4}$ & 2 & 1 \\ \hline
			4.02 & \vspace{0.1 mm} $ -A^6 K_1 - A^4 K_1^2 +2 K_2  +3 -2K_1 ^2  + A^2 K_1  + A^{-2} K_1 - A^{-4} K_1^2 - A^{-6} K_1 $ & 2 & 1 \\ \hline
			4.03 & \vspace{0.1 mm} $ A^{8}K_1^{2} -A^{4} - K_1^{2} +1 -2A^{2} K_1 ^{1} - A^{4} K_2 + K_2 +2A^{-2} K_1^{1} + A^{-4}$ & 2 & 1 \\ \hline
			4.04 & \vspace{0.1 mm} $A^{2} - A^{4}K_1^{1} -2A^{2}K_1^{2} -2A^{-2}K_1^{2} + A^{-2}K_2 +A^{2} K_2 +2A^{-2} +1K_1$ & 2 & 1 \\ \hline
			4.05 & \vspace{0.1 mm} $ -A^{4}K_1^{1} + A^{2} -2A^{2}K_1^{2} -2A^{-2}K_1^{2} +A^{-2}K_2 +A^{2}K_2 + K_1 +2A^{-2} $ & 2 & 1 \\ \hline
			4.06 & \vspace{0.1 mm} $ -A^{6}K_1 -A^{4}K_1^{2} +K_2 + A^{-4}K_1^{2} 
      -A^{-4}K_2 + A^{2}K_1 +2 + A^{-2}K_1  - A^{-4} - A^{-6}K_1  $ & 2 & 1 \\ \hline
			4.07 & \vspace{0.1 mm} $ A^{8}K_1^{2} -3K_1^{2} +2K_2^{1} -2A^{4}K_1^{2} +2 -2A^{2}K_1^{1} +2A^{-2}K_1^{1} +A^{-4} $ & 2 & 1 \\ \hline
		  4.08 & \vspace{0.1 mm} $ -A^{6}K_1 -A^{4}K_1^{2} +3 +2K_2 -2K_1^{2} +A^{2}K_1 +A^{-2}K_1 -A^{-4}K_1^{2} -A^{-6}K_1 $ & 2 & 1 \\ \hline
		  4.09 & \vspace{0.1 mm} $ -A^{4} -A^{2}K_1 +A^{-2}K_1 +A^{-4} $ & 1 & 1 \\ \hline
			4.10 & \vspace{0.1 mm} $ -A^{6} -A^{4}K_1 +2A^{2} +2K_1 +A^{-2}
      -A^{2}K_1^{2} -A^{-2}K_1^{2} +A^{-2}K_2 -A^{-4}K_1 $ & 2 & 1 \\ \hline
			4.11 & \vspace{0.1 mm} $ -A^{-2}K_1^{2} +A^{2} -A^{4}K_1 -A^{2}K_1^{2} +A^{-2}K_2 + K_1 +A^{-2} $ & 2 & 1 \\ \hline
			4.12 & \vspace{0.1 mm} $ -A^{6}K_1 -A^{4}K_2 +A^{2}K_1 +1 +2K_2 -A^{-4}K_2 +A^{-2}K_1 -A^{-6}K_1 $ & 2 & 1 \\ \hline
			4.13 & \vspace{0.1 mm} $ A^{4} -A^{-4}K_1^{2} +1 -A^{4}K_1^{2} -2K_1^{2} +2K_2 +A^{-4} $ & 2 & 1 \\ \hline
			4.14 & \vspace{0.1 mm} $ -A^{6}K_1 -2A^{-4}K_1^{2} +A^{2}K_1 +2K_2 +2 -2K_1^{2} +A^{-4} $ & 2 & 1 \\ \hline
			4.15 & \vspace{0.1 mm} $ K_2 -A^{4}K_1^{2} -K_1^{2} +1 -A^{2}K_1 +A^{-2}K_1 +A^{-4}$ & 2 & 1 \\ \hline
		\end{tabular}
	\caption{Bounds: Knots 1-15 }
	\label{tab:Bounds0}
\end{table*}

\begin{table*}
		\begin{tabular}{l|p{3 in}|l|l}
			Knot & Arrow Polynomial & v(K) & g(K) \\ \hline
			4.16 & \vspace{0.1 mm} $ -A^{6} -A^{4}K_1 +A^{2} +2K_1 +A^{-2} -A^{-4}K_1 $ & 1 & 1 \\ \hline
			4.17 & \vspace{0.1 mm} $ -A^{6}K_1 -A^{4}K_2 +A^{2}K_1 +2 +2K_2 -A^{-4}K_1^{2} -K_1^{2} +A^{-2}K_1 -A^{-6}K_1 $ & 2 & 1 \\ \hline
			4.18 & \vspace{0.1 mm} $ -2A^{-2}K_1^{2} -A^{4}K_1 +A^{2}K_2 -2A^{2}K_1^{2} +2A^{-2} +A^{2} +K_1 +A^{-2}K_2 $ & 2 & 1 \\ \hline
			4.19 & \vspace{0.1 mm} $ A^{4} + K_2 +1 -A^{4}K_1^{2} -K_1^{2} $ & 2 & 1 \\ \hline
			4.20 & \vspace{0.1 mm} $ -A^{6}K_1 +K_2 +A^{2}K_1 +2 -A^{-4}K_1^{2} -K_1^{2} $ & 2 & 1 \\ \hline
			4.21 & \vspace{0.1 mm} $ A^{2} +A^{2}K_2 -A^{4}K_1 
      -A^{2}K_1^{2} -2A^{-2}K_1^{2} +2A^{-2} +2K_1 +A^{-2}K_2 -A^{-4}K_1 -A^{-6}K_1^{2}$ & 2 & 1 \\ \hline
       4.22
& \vspace{0.1 mm} $-A^{6}K_1 +A^{2}K_1 +K_2 +2 -A^{2}K_1K_2 -A^{-2}K_1K_2 
      -K_1^{2} -A^{-4}K_1^{2} +A^{2}K_1^{3} +2A^{-2}K_1^{3} +A^{-6}K_1^{3} -A^{-2}K_1 -A^{-6}K_1$ & 3 & 2 \\ \hline
			4.23 & \vspace{0.1 mm} $A^{8}K_1 -2A^{4}K_1 +2A^{-2} -A^{2}K_1^{2} -A^{-2}K_1^{2} +K_1 +A^{-2}K_2 $& 2 & 1 \\ \hline
			4.24
& \vspace{0.1 mm} $ A^{2} -2A^{4}K_1 +A^{4}K_1^{3} +2K_1^{3} +A^{-4}K_1^{3} -A^{-4}K_1 
      +A^{-2} + K_1 +A^{-2}K_2 -1K_1K_2 -A^{-4}K_1K_2 -A^{-2}K_1^{2} -A^{-6}K_1^{2}$ &3&2 \\ \hline
 4.25 & \vspace{0.1 mm} $-A^{6}K_1 -A^{-2}K_1 +A^{-6}K_1 +A^{2}K_1 -A^{-4} +1 +A^{-8} $& 1 & 1\\ \hline
	4.26 & \vspace{0.1 mm} $ A^{2}K_1 -A^{-2}K_1K_2 -A^{2}K_1K_2 +1 +A^{-2}K_3$ & 3 & 2 \\ \hline
4.27 & \vspace{0.1 mm} $A^{-6} -A^{-2} +A^{2} -A^{-4}K_1 +A^{-8}K_1$ & 1& 1\\ \hline
4.28 & \vspace{0.1 mm} $ A^{2} +K_3  + K_1 -A^{-4}K_1K_2 - K_1K_2  $&3&2\\ \hline
4.29 & \vspace{0.1 mm} $1 -A^{4}K_1^{2} -K_1^{2} +K_2 -A^{2}K_1 +A^{-2}K_1 +A^{-4} $ & 2 & 1 \\ \hline
 4.30 & \vspace{0.1 mm} $-2A^{-2}K_1^{2} +A^{2}K_2 +A^{2} -2A^{2}K_1^{2} +A^{-2}K_2 -A^{4}K_1 +2A^{-2} + K_1$ & 2 & 1 \\ \hline
		\end{tabular}
	\caption{Bounds: Knots 16-30 }
	\label{tab:Bounds1}
\end{table*}

\begin{table*}
		\begin{tabular}{l|p{3 in}|l|l}
			Knot & Arrow Polynomial & v(K) & g(K) \\ \hline
			4.31 & \vspace{0.1 mm} $ -A^{6}K_1^{2} -2A^{2}K_1^{2} +2A^{2}K_2 -A^{-2}K_1^{2} 
      -A^{4}K_1 +2K_1 +2A^{-2} +A^{2} -A^{-4}K_1 $ & 2 & 1 \\ \hline
4.32 & \vspace{0.1 mm} $-A^{6}K_1 +A^{2}K_1 +3 -A^{4}K_1^{2} - K_1^{2} -A^{-4} +K_2 +A^{-2}K_1 -A^{-6}K_1 $ & 2 & 1 \\ \hline
4.33 & \vspace{0.1 mm} $ -A^{4}K_1 + K_1 +A^{-2}$ & 1 & 1 \\ \hline
4.34 & \vspace{0.1 mm} $ -A^{6}K_1 +A^{2}K_1 +2 -A^{-4}K_1^{2} +K_2 - K_1^{2} $& 2 & 1 \\ \hline
4.35 & \vspace{0.1 mm} $ - K_1^{2} -A^{-4}K_1^{2} +2 +A^{-4}K_2$ & 2 & 1 \\ \hline
4.36 & \vspace{0.1 mm} $ A^{2} +2K_1 -A^{4}K_1 -A^{-4}K_1 +A^{-2}K_2 -A^{-6}K_2$ & 2 & 1 \\ \hline
4.37 & \vspace{0.1 mm} $K_2 -A^{4}K_2 -A^{2}K_1 +A^{-2}K_1 +A^{-4}$ & 2 & 1 \\ \hline
4.38 & \vspace{0.1 mm} $A^{2}K_2 -A^{4}K_1 -A^{2}K_1^{2} -A^{-2}K_1^{2} +2A^{-2} + K_1 $& 2 & 1 \\ \hline
4.39 & \vspace{0.1 mm} $A^{2}K_2 -A^{4}K_1K_2 -K_1K_2 -A^{2}K_1^{2} -A^{-2}K_1^{2} 
      -A^{4}K_1 +2A^{-2} -A^{-4}K_1 +A^{4}K_1^{3} +2K_1^{3} +A^{-4}K_1^{3} $& 3 & 2 \\ \hline
      4.40 & \vspace{0.1 mm} $-A^{6}K_1 -A^{-4} +2 +A^{2}K_1 $& 1 & 1 \\ \hline
		4.41 & \vspace{0.1 mm} $A^{8}K_1 -2A^{4}K_1 -A^{2} + K_1 +2A^{-2}$ & 1 & 1 \\ \hline
4.42 & \vspace{0.1 mm} $
 -A^{6}K_1 +A^{6}K_1^{3} +2A^{2}K_1^{3} +A^{-2}K_1^{3} -A^{2}K_1 +2 -A^{2}K_1K_2 -A^{-2}K_1K_2 - K_1^{2} -A^{-4}K_1^{2} +A^{-4}K_2$ & 2 & 2 \\ \hline 
4.43
& \vspace{0.1 mm} $ -A^{6}K_1 -A^{-2}K_1 +A^{-6}K_1 +A^{2}K_1 -A^{-4} +1 +A^{-8}$
& 1 & 1 
\\ \hline 4.44
& \vspace{0.1 mm} $-A^{4}K_1 +A^{-2} + K_1 $
& 1 & 1
\\ \hline 4.45
& \vspace{0.1 mm} $ -K_1K_2 -A^{4}K_1K_2 + K_1 +K_3 +A^{-2}$
3 & 1
\\ \hline 
\end{tabular}
	\caption{Bounds: Knots 30-45 }
	\label{tab:Bounds2}
\end{table*}

\begin{table*}
		\begin{tabular}{l|p{3 in}|l|l}
			Knot & Arrow Polynomial & v(K) & g(K) \\ \hline
4.46
& \vspace{0.1 mm} $1$ & 0 & 0 
\\ \hline 4.47
& \vspace{0.1 mm} $ A^{2}K_3^{1} +1 -A^{-2}K_1K_2 -A^{2}K_1K_2 +A^{-2}K_1$
& 3 & 2 \\ \hline
4.48
& \vspace{0.1 mm} $A^{4} -2A^{4}K_1^{2} -2K_1^{2} +1 +2K_2 -A^{2}K_1 +A^{-2}K_1 +A^{-4}$
& 2 & 1 
\\ \hline 4.49
& \vspace{0.1 mm} $A^{2}K_2 -A^{-2}K_1^{2} -A^{2}K_1^{2} -A^{4}K_1 + K_1 +2A^{-2} $
& 2 & 1 
\\ \hline 4.50
& \vspace{0.1 mm} $-A^{6}K_1^{2} -A^{2}K_1^{2} +A^{2} +2K_1 +A^{2}K_2 -A^{4}K_1 -A^{-4}K_1 +A^{-2}$
& 2 & 1 
\\ \hline 4.51
& \vspace{0.1 mm} $-A^{6}K_1 +A^{2}K_1 +3 -A^{4}K_1^{2} -2K_1^{2} +2K_2 -A^{-4}K_1^{2} +A^{-2}K_1 -A^{-6}K_1$
& 2 & 1 
\\ \hline 4.52
& \vspace{0.1 mm} $ -A^{6}K_1 +A^{2}K_1 +2 -A^{-4}$
& 1 & 1 \\ \hline
4.53
& \vspace{0.1 mm} $ A^{8} -2A^{4} -2A^{2}K_1 +1 +2A^{-2}K_1 +A^{-4}$
& 1 & 1 
\\ \hline 4.54
& \vspace{0.1 mm} $ -A^{6} -A^{4}K_1 +A^{6}K_1^{2} +A^{-2} -A^{2}K_2 +A^{2} + K_1 -A^{-2}K_1^{2} +A^{-2}K_2$
& 2 & 1 
\\ \hline 4.55
& \vspace{0.1 mm} $ A^{4} +2K_2 +1 -A^{4}K_1^{2} -2K_1^{2} -A^{-4}K_1^{2} +A^{-4}$
& 2 & 1 
\\ \hline 4.56
& \vspace{0.1 mm} $ A^{4} +1 -A^{4}K_1^{2} +2K_2 -2K_1^{2} -A^{-4}K_1^{2} +A^{-4}$
& 2 & 1 
\\ \hline 4.57
& \vspace{0.1 mm} $-A^{6}K_2 -A^{4}K_1 +2A^{2}K_2 +2K_1 +2A^{-2} -A^{2}K_1^{2} -A^{-2}K_1^{2} -A^{-4}K_1$
& 2 & 1 
\\ \hline 4.58
& \vspace{0.1 mm} $ -A^{6}K_1 +A^{2}K_1 +3 -A^{4} +A^{-2}K_1 -A^{-4} -A^{-6}K_1$
& 1 & 1 
\\ \hline 4.59
& \vspace{0.1 mm} $ A^{4}K_2 -A^{-4}K_1^{2} -A^{4}K_1^{2} -2K_1^{2} +3 +A^{-4}K_2$
& 2 & 1 
\\ \hline 4.60
& \vspace{0.1 mm} $ -A^{6}K_1 -2A^{-4}K_1^{2} +A^{2}K_1 +3 + K_2 -2K_1^{2} +A^{-4}K_2$
& 2 & 1 \\ \hline
		\end{tabular}
	\caption{Bounds: Knots 46-60 }
	\label{tab:Bounds3}
\end{table*}

\begin{table*}
		\begin{tabular}{l|p{3 in} | l | l }
			Knot & Arrow Polynomial & v(K) & g(K) \\ \hline
 4.61
& \vspace{0.1 mm} $1 -A^{4} -A^{2}K_1 +A^{-2}K_1 +A^{-4}$ & 1 & 1 
\\ \hline 4.62
& \vspace{0.1 mm} $A^{2} -A^{4}K_1K_2 -K_1K_2 -A^{2}K_1^{2} -A^{-2}K_1^{2} 
      -A^{4}K_1 +A^{-2}K_2 -A^{-4}K_1 +A^{-2} +A^{4}K_1^{3} +2K_1^{3} +A^{-4}K_1^{3}$
& 2 & 1 
\\ \hline 4.63
& \vspace{0.1 mm} $ A^{2} -A^{4}K_1 -A^{2}K_1^{2} -A^{-2}K_1^{2} +A^{-2}K_2 +K_1 +A^{-2}$
& 2 & 1
\\ \hline 4.64
& \vspace{0.1 mm} $ -A^{6}K_1 -A^{-4}K_2 +K_2 +1 +A^{2}K_1$
&2 & 1 \\ \hline
 4.65
& \vspace{0.1 mm} $ A^{8}K_1 -2A^{4}K_1 +A^{-2} -A^{2}K_2 +K_1 +A^{-2}K_2$
& 2 & 1 \\ \hline
4.66
& \vspace{0.1 mm} $-A^{6}K_1 +A^{-2}K_1^{3} +1 +A^{6}K_1^{3} +2A^{2}K_1^{3} -A^{-4}K_1^{2}
      +K_2 -A^{2}K_1 -A^{2}K_1K_2 -A^{-2}K_1K_2 -K_1^{2} +A^{-4}$
& 3 & 2 \\ \hline
4.67
& \vspace{0.1 mm} $ -A^{-4}K_1^{2} +1 +K_2 -K_1^{2} +A^{-4}$
& 2 & 1 
\\ \hline 4.68
& \vspace{0.1 mm} $ A^{2} -A^{4}K_1 +A^{-2} +2K_1 -A^{-4}K_1 -A^{-6}$
& 1 & 1 
\\ \hline 4.69
& \vspace{0.1 mm} $ A^{4}K_2 -2A^{4}K_1^{2} -2K_1^{2} +K_2 +2 -A^{2}K_1 +A^{-2}K_1 +A^{-4}$
& 2 & 1 
\\ \hline 4.70
& \vspace{0.1 mm} $-A^{6}K_1^{2} -A^{2}K_1^{2} +A^{2}K_2 +2K_1 -A^{4}K_1 +A^{2} -A^{-4}K_1 +A^{-2}$
& 2 & 1 
\\ \hline 4.71
& \vspace{0.1 mm} $ -A^{6}K_1 +A^{2}K_1 +2K_2 -A^{4}K_1^{2} -2K_1^{2} +3 -A^{-4}K_1^{2} +A^{-2}K_1 -A^{-6}K_1$
2 & 1 
\\ \hline 4.72
& \vspace{0.1 mm} $1$
& 0 & 0 
\\ \hline 4.73
& \vspace{0.1 mm} $A^{8}K_2 -2A^{4}K_2 -2A^{2}K_1 +K_2 +2A^{-2}K_1 +A^{-4}$ &
2 & 1
\\ \hline 4.74
& \vspace{0.1 mm} $-A^{6}K_2 -A^{4}K_1 +A^{6}K_1^{2} +2A^{-2} -A^{2} +A^{2}K_2 +K_1 -A^{-2}K_1^{2}$ &
2 & 1 \\ \hline
4.75
& \vspace{0.1 mm} $ -A^{6}K_1 -A^{4} +3 +A^{2}K_1 +A^{-2}K_1 -A^{-4} -A^{-6}K_1$
& 1 & 1 
\\ \hline 
		\end{tabular}
	\caption{Bounds: Knots 61-75 }
	\label{tab:Bounds4}
\end{table*}

\begin{table*}
		\begin{tabular}{l|p{3 in} | l | l }
			Knot & Arrow Polynomial & v(K) & g(K) \\ \hline
		 4.76
& \vspace{0.1 mm} $ A^{4}K_2 -A^{4}K_1^{2} +3 -2K_1^{2} -A^{-4}K_1^{2} +A^{-4}K_2$
& 2 & 1 
\\ \hline 4.77
& \vspace{0.1 mm} $ A^{4}K_2 -2K_1^{2} -A^{4}K_1^{2} +3 -A^{-4}K_1^{2} +A^{-4}K_2$
& 2 & 1 
\\ \hline 4.78
& \vspace{0.1 mm} $ -A^{6}K_1K_2 -A^{2}K_1K_2 -A^{4}K_1^{2} -K_1^{2} +A^{-2}K_1^{3} 
      +K_2 +1 +A^{6}K_1^{3} +2A^{2}K_1^{3} -2A^{2}K_1 +A^{-4}$
& 3 & 2 
\\ \hline 4.79
& \vspace{0.1 mm} $ A^{8}K_1^{3} +2A^{4}K_1^{3} +K_1^{3} -3A^{4}K_1 
      +2A^{-2} -A^{4}K_1K_2 -K_1K_2 -A^{2}K_1^{2} -A^{-2}K_1^{2} +K_1 +A^{-2}K_2$
& 3 &2 
\\ \hline 4.80
& \vspace{0.1 mm} $-A^{6}K_1K_2 -A^{2}K_1K_2 -A^{-2}K_1 
      +A^{2}K_1 +A^{-6}K_1 +A^{2}K_3 -A^{-4} +1 +A^{-8}$
&3 & 2 \\ \hline
4.81
& \vspace{0.1 mm} $ A^{4}K_3 -A^{-2} +A^{2} +A^{-6} -A^{4}K_1 K_2 -K_1 K_2 -A^{-4}K_1 +K_1 +A^{-8}K_1$
& 3 & 2 
\\ \hline 4.82
& \vspace{0.1 mm} $ -A^{6}K_1 +A^{-4}K_1^{2} -A^{4}K_1^{2} +2 -2A^{-4} +A^{2}K_1 +A^{-8}$
& 2 & 1 
\\ \hline 4.83
& \vspace{0.1 mm} $ -K_1 K_2 -A^{4}K_1 K_2 +K_1 +A^{-2} +K_3$
& 3 & 2
\\ \hline 4.84
& \vspace{0.1 mm} $ 2 +A^{2}K_1 -A^{-4} -A^{-2}K_1 -K_1^{2} +A^{-8}K_1^{2}$
& 2 & 1
\\ \hline 4.85
& \vspace{0.1 mm} $ -A^{2}K_1^{2} +A^{2} +A^{-6}K_1^{2}$
& 2 & 1 
\\ \hline 4.86
& \vspace{0.1 mm} $ A^{8} -A^{4} +2 -A^{-4} -K_1^{2} +A^{-8}K_1^{2}$
& 2 & 1
\\ \hline 4.87
& \vspace{0.1 mm} $ -A^{6}K_1 K_2 -A^{2}K_1 K_2 -A^{4}K_1^{2} -2A^{-4} 
       +A^{2}K_1 +2 +A^{2}K_3 +A^{-4}K_1^{2} +A^{-8}$
& 3 & 2 
\\ \hline 4.88
& \vspace{0.1 mm} $ 2 +A^{2}K_1 +A^{2}K_3 -A^{-4} -A^{2}K_1 K_2 
      -A^{-2}K_1 K_2 -K_1^{2} +A^{-8}K_1^{2}$
& 3 & 2 
\\ \hline 4.89
& \vspace{0.1 mm} $ A^{4} -2A^{4}K_1^{2} +2A^{-4}K_1^{2} -2A^{-4} +1 +A^{-8}$
& 2 & 1
\\ \hline 4.90
& \vspace{0.1 mm} $ A^{8}K_1^{2} -A^{4} -2K_1^{2} +3 -A^{-4} +A^{-8}K_1^{2}$
& 2 & 1
\\ \hline 
		\end{tabular}
	\caption{Bounds: Knots 76-90 }
	\label{tab:Bounds5}
\end{table*}

\begin{table*}
		\begin{tabular}{l|p{3 in} | l | l }
			Knot & Arrow Polynomial & v(K) & g(K) \\ \hline
	4.91
& \vspace{0.1 mm} $ -A^{10}K_1^{3} -A^{6}K_1^{3} +A^{2}K_1^{3} +2A^{6}K_1 
       +A^{-2}K_1^{3} -2A^{2}K_1 +A^{-4}$
& 3 & 1
\\ \hline 4.92
& \vspace{0.1 mm} $-A^{6}K_3 -A^{4} +2 +A^{2}K_3 -A^{-4} +A^{-8}$
& 3 & 1
\\ \hline 4.93
& \vspace{0.1 mm} $ A^{4}K_3 +A^{-6}K_1^{2} +A^{2} -A^{2}K_1^{2} -A^{4}K_1 K_2 -K_1 K_2 +K_1$
& 3 & 2
\\ \hline 4.94
& \vspace{0.1 mm} $ -A^{6}K_1 +2 -A^{4} +A^{2}K_1 -A^{-4} +A^{-8}$
& 1 & 1
\\ \hline 4.95
& \vspace{0.1 mm} $ -A^{4}K_3 +A^{-2} +K_3$
& 3 & 1
\\ \hline 4.96
& \vspace{0.1 mm} $ A^{-6}K_1^{2} -A^{2}K_1^{2} +A^{2}$
& 2 & 1
\\ \hline 4.97
& \vspace{0.1 mm} $ -A^{2}K_1 K_2 -A^{-2}K_1 K_2 +A^{2}K_3 +1 +A^{-2}K_1$
& 3 & 1 
\\ \hline 4.98
& \vspace{0.1 mm} $ 1$
& 0 & 0
\\ \hline 4.99
& \vspace{0.1 mm} $ A^{8} -A^{4} -A^{-4} +1 +A^{-8}$
& 0 & 0
\\ \hline 4.100
& \vspace{0.1 mm} $ -A^{10}K_1 +A^{6}K_1 -A^{2}K_1 +A^{-2}K_1 +A^{-4}$
& 1 & 1 \\ \hline
	4.101
& \vspace{0.1 mm} $ A^{8}K_1 +K_1^{3} -A^{-4}K_1 +A^{-2} -A^{8}K_1^{3} -A^{4}K_1^{3} +A^{-4}K_1^{3}$
& 3 & 1 
\\ \hline 4.102
& \vspace{0.1 mm} $ -A^{6}K_1 -A^{-2}K_1^{3} +1 +A^{6}K_1^{3} +A^{2}K_1^{3}
      -A^{-6}K_1^{3} -A^{2}K_1 +A^{-6}K_1 +A^{-2}K_1$
& 3 & 1
\\ \hline 4.103
& \vspace{0.1 mm} $A^{4}K_1 +A^{2} +K_3 -A^{4}K_1 K_2 -K_1 K_2 -A^{2}K_1^{2} +A^{-6}K_1^{2}$
& 3 & 2
\\ \hline 4.104
& \vspace{0.1 mm} $ A^{2}K_3 -A^{-2}K_3 -A^{-4} +1 +A^{-8}$
& 3 & 1
\\ \hline 4.105
& \vspace{0.1 mm} $-A^{4} +1 +A^{-8}$
& 0 & 0 
\\ \hline 4.106
& \vspace{0.1 mm} $ A^{2} -A^{2}K_1^{2} +A^{-6}K_1^{2}$
& 2 & 1 
\\ \hline 4.107
& \vspace{0.1 mm} $1$
& 0 & 0 
\\ \hline 4.108
& \vspace{0.1 mm} $A^{8} -A^{4} -A^{-4} +1 +A^{-8}$
& 0 & 0 \\ \hline
		\end{tabular}
	\caption{Bounds: Knots 90-108 }
	\label{tab:Bounds6}
\end{table*}

\end{document}